\documentclass[12pt,twoside]{article}
\usepackage{amsmath,amsfonts,amssymb,amsthm,amscd}
\makeindex

\usepackage{amssymb}
\usepackage[latin1]{inputenc}
\usepackage{ epsfig,epsf}
\usepackage{enumerate}
\usepackage{indentfirst}
\usepackage{euscript}
\usepackage{graphicx}
\usepackage{amsmath}
\usepackage{amsfonts}
\usepackage{amssymb}
\makeindex

 \theoremstyle{plain}

\newtheorem{lemma}{Lemma}

\theoremstyle{definition}
\newtheorem{definition}{Definition}
\newtheorem{theorem}{Theorem}

\newcommand{\ddd}{\displaystyle}
\newcommand{\nd}{\noindent}
\newcommand{\vsp}{\vspace{0.2cm}}
\newcommand{\fimaf}{$\hfill{\Box}$}
\newcommand{\fim}{$\hfill{\rule{2.5mm}{2.5mm}}$}
\usepackage{latexsym}
\theoremstyle{plain}

\begin{document}


{\large \bf On the study of solutions for a non linear differential equation on compact Riemannian Manifolds.}


\begin{center}{{C.R. Silva}
 and  {Marcelo Souza}}
\end{center}

\begin{abstract}
 { In this paper we study the existence of solutions for a class of non-linear differential  equation on compact Riemannian manifolds. 
 We establish a lower and upper solutions' method to show the existence of a smooth positive solution for the equation (\ref{E4})
\begin{equation}
 \label{E4}
\Delta u \ + \ a(x)u \ = \ f(x)F(u) \ + \ h(x)H(u),
\end{equation}

\nd where \ $a, \ f, \ h$ \  are positive smooth functions on $M^n$, a $n-$dimensional compact Riemannian manifold,
 and \ $ F, \ H$ 
\ are non-decreasing smooth functions on $\mathbb{R}$.  In \cite{djadli} the equation (\ref{E4}) was studied when
$F(u)=u^{2^{\ast}-1} $ and $H(u)=u^q$ in the Riemannian context, i.e.,
\begin{equation}
\label{E3}
\Delta u \ + \ a(x)u \ = \ f(x)u^{2^{\ast}-1} \ + \ h(x)u^q, 
\end{equation}
\nd where \ $0 \ < \ q \ < 1$. In \cite{correa} Corr\^ea,  Gon\c{c}alves and Melo  studied
an equation of the type equation (\ref{E3}),  in the Euclidean context.

}\end{abstract}

\thispagestyle{empty}
{ \small Subjclass: [2010]{Primary 53C21; Secondary 35J60}



{The authors were partially supported by CAPES}}

\section{ Introduction}

The study of the theory of nonlinear differential equations on 
Riemannian manifolds, has began in 1960 with the so-called Yamabe problem. At a time when little was known 
about the methods of attacking a non-linear equation, the Yamabe problem came to light of a geometric idea and from time sealed a 
merger of the areas of geometry and differential equations.
Let \ $(M,g)$ be a compact Riemannian manifold of dimension \ $n$, \ $n \geq 3$. \ Given \ $\widetilde{g} = u^{4/(n - 2)}g$
 \  some conformal metrical to the metric \ $g$, \ is well known that the scalar curvatures
  \ $R$ \ and \ $\widetilde{R}$ \ of the metrics \ $g$ \ and \ $\widetilde{g}$, 
respectively, satisfy the law of transformation
$$\Delta u \ + \ \ddd\frac{n - 2}{4(n - 1)}Ru \ = \ \ddd\frac{n - 2}{4(n-1)}\widetilde{R}u^{2^{\ast}-1}$$

\nd where \ $\Delta$ \ denote the Laplacian operator associated to \ $g$.

In 1960, Yamabe \cite{yamabe} announced that for every compact Riemannian manifold \ $(M,g)$ 
\ there exist  a metric \ $\widetilde{g}$ \ conformal to \ $g$ \ for which \ $\widetilde{R}$ \ is constant.
 In another words, this mean that for every compact Riemannian manifold \ $(M,g)$ there exist 
\ $u \in C^{\infty}(M), \ u > 0 \ \mbox{on} \ M$ \ and \ $\lambda \in \mathbb{R}$ \ such that  
$$\Delta u \ + \ \ddd\frac{n - 2}{4(n - 1)}Ru \ = \ \lambda u^{2^{\ast}-1}. \eqno {(Y)}$$

In 1968, Tr$\ddot{u}$dinger \cite{trudinger} found an error in the work of Yamabe, which generated a race to
solve what became known as the Yamabe problem, today it is completely positively resolved,
that is, the assertion of  Yamabe is true.

The main step towards the resolution of the Yamabe problem was given in 1976 by Aubin in his classic article \cite{aubin2}.
 In \cite{aubin2} Aubin showed that the statement was true since the manifold satisfy a condition
on an invariant (called Yamabe invariant).
Then he used  tests functions, locally defined to show that  non locally conformal flat manifolds, of
dimension  $ n > 6$, \ satisfying this condition. Finally, the problem for $n\geq 3 $ was completed  solved by R. Schoen \cite{schoen}. 

As previously reported, several disturbances were the Yamabe problem, all of analytical character, both in the sense of equation 
(with the addition of other factors) and in the sense of the operator (the Laplacian for the $ p $-Laplacian), and all (at least
 those in this study) using the idea of estimating the functional Aubin corresponding multiple functions \ $ u_ {\lambda} $. We
 can cite some articles, such as \cite {azorero1}, \cite {brezis1}, \cite {demegel}, \cite {djadli}, \cite {druet4} and 
\cite {olimpio} .

This work aims to work with problems related to the equation 
$ (Y) $, although, as we shall see, with different methods from those used by Yamabe, these results were obtained in \cite{silva}, and some of them were published in \cite{silva2}. 
The above equation was studied  by Gon\c{c}alves and Alves \cite{correa}, in the  Euclidean space,
 here  by using the methods in \cite{correa} we study the case  in the Riemannian space.

\begin{definition}
We  say that a function \ $\underline{u}$ \ (respectively, \ $\overline{u}$) \ $ \in \ H^2_1\cap L^{\infty}$, \linebreak
$\underline{u} \ \geq \ 0 \ (\overline{u} \ \geq \ 0)$ \ is a lower solution (respectively, upper solution)  of equation (\ref{E4}) 
if for all \ $\varphi \ \in \ H^2_1, \ \varphi \ \geq \ 0$
$$\ddd\int_M{\nabla \underline{u}.\nabla \varphi dV} \ + \ \ddd\int_M{a \underline{u} \varphi dV} \ \leq \ \ddd\int_M{f(x) 
F(\underline{u}) \varphi dV} \ + \ \ddd\int_M{h(x) H(\underline{u}) \varphi dV}$$ 

\nd (respectively,
$$\ddd\int_M{\nabla \overline{u}.\nabla \varphi dV} \ + \ \ddd\int_M{a \overline{u} \varphi dV} \ \geq \ \ddd\int_M{f(x) 
F(\overline{u}) \varphi dV} \ + \ \ddd\int_M{h(x) H(\overline{u}) \varphi dV}).$$ 
\end{definition}

We consider the conditions:

\vspace{0.35cm}
\nd $(\alpha_1)  \left\{ 
\begin{array}{lcl}
0 \leq F(t) \leq t^{2^{\ast}-1} \ \ \mbox{and} \ \ 0  \leq  H(t)  \leq  t^q  \ \ & \mbox{if} & t \geq  0 \\
F(t) = H(t) = 0 & \mbox{if} & t < 0
\end{array}
\right. $
\vsp

\nd where \ \ $2^{\ast}  =  2n/(n-2) \ \ \mbox{and} \ \ 0  <  q  < 2^{\ast} -1.$

\vspace{0.35cm}

\nd $(\alpha_2) \  a  >  0, \ f  \geq  0, \ f  \not\equiv  0 \ \ \mbox{and} \ \ h  \geq  0, \ \ \ h  \not\equiv  0$.

We proved, in the main result
\begin{theorem}
Let \ $(M,g)$ be a compact $n-$dimensional Riemannian manifold $n\geq 3$. Suppose that \ $(\alpha_1)$ \ and \ $(\alpha_2)$ holds.
 If \ $\underline{u}, \ \overline{u} \ \in \ H^2_1\cap L^{\infty}$ \ 
are, respectively, the lower and the upper solutions  of the equation \ (\ref{E4}) \ with \ $0 \ \leq 
\ \underline{u} \ \leq \ \overline{u}$
 \ a.e. in \ $M$ \ and \ $\underline{u} \ \not\equiv \ 0$, \ then the equation \ (\ref{E4})
$$\Delta u \ + \ a(x)u \ = \ f(x)F(u) \ + \ h(x)H(u)$$ \ admits a  positive 
regular solution  \ $u$, \ such that $\underline{u} \ \leq \ u \ \leq \ \overline{u}$.
\vsp

\vsp\end{theorem}

\section{An auxiliary Lemma}
  
We will need an auxiliary lemma.

Let us consider the  following equation: 
\begin{equation}
\label{4.1}
\Delta u \ + \ a(x)u \ = \ \Psi \ \ \mbox{ in } \ M. 
 \end{equation}

\begin{lemma}
Assume that \ $a \ > \ 0$. If \ $\Psi \ \in \ (H^2_1)^{'}$, \ then the equation \ (\ref{4.1}) has a unique solution
\ $u \ \in \ H^2_1$.  Moreover, the $T$ operator 

\begin{eqnarray*}
T \ : & (H^2_1)^{'} & \longrightarrow \ H^2_1 \\
      &  \Psi      &  \longrightarrow \ T(\Psi) \ = \ u
\end{eqnarray*}

\nd is continuous and non-decreasing.
\end{lemma}

\nd {\bf Proof of the Lemma 1:} 

Consider the functional \ $I \ : \ H_1^2 \ \longrightarrow \ \mathbb{R}$ \ given by
$$I(u) \ = \ \ddd\frac{1}{2}\int_M{|\nabla u|^2dV} \ + \ \ddd\frac{1}{2}\int_M{au^2dV} \ - \ \bigl<\Psi, u\bigr> , \ \ 
\ u \ \in \ H^2_1.$$

Then, \ $I \ \in \ C^1(H^2_1,\mathbb{R})$ \ \ and  
$$\bigl<I^{'}(u),v\bigr> \ = \ \ddd\int_M{\nabla u\nabla vdV} \ + \ \ddd\int_M{auvdV} \ - \ \bigl<\Psi, v\bigr> , \ \ 
\ u, \ v \ \in \ H^2_1.$$

Therefore, the  solutions of the equation \ (\ref{4.1}) are critical points of the functional \ $I$.

As \ $a \ > \ 0$ \ and \ $M$ \ is a compact manifold, 
$$I(u) \ \geq \ \ddd\frac{1}{2}\int_M{|\nabla u|^2dV} \ + \ \ddd\frac{\min a}{2}\int_M{u^2dV} \ - \ \bigl<\Psi, u\bigr> 
\ \geq \ C\|u\|^2_{H^2_1} \ - \ \bigl<\Psi,u\bigr> \ , \ u \ \in \ H^2_1,$$

\nd where \ $C \ = \ \min{\{1/2,(\min a)/2\}} \ > \ 0$.

Since \ $\bigl<\Psi,u\bigr> \ \leq \ \|\Psi\|_{(H^2_1)^{'}}\|u\|_{H^2_1}$, \ we have that
$$I(u) \ \geq \ C\|u\|^2_{H^2_1} \ - \ \|\Psi\|_{(H^2_1)^{'}}\|u\|_{H^2_1} \ \ \forall \ u \ \in \ H^2_1.$$

Then, \ $I$ \ is coercive.

{{\bf Claim 1}
\it The functional $I$ \ is l.s.s.c. (lower sequentially semi continuous), namely, for all \ $u_i \ \rightharpoonup \ u$
 \ in \ $H^2_1$ \ implies that \ $I(u) \ \leq \ \ddd\lim_{i \to \infty}{\inf{I(u_i)}}$.
}

Indeed, suppose that \ $u_i \ \rightharpoonup \ u$ \ in \ $H^2_1$.

Since the embedded \ $H^2_1 \ \hookrightarrow \ L^2$ \ is compact and \ $H^2_1$ \ is reflexive, 
$$u_i \ \longrightarrow \ u \ \mbox{in} \ L^2,$$

\nd and

$$\|u\|^2_{H^2_1} \ \leq \ \ddd\lim_{i \to \infty}{\inf{\|u_i\|^2_{H^2_1}}}.$$

With this,
$$\ddd\int_M{|\nabla u|^2dV} \ \leq \ \ddd\lim_{i \to \infty}{\inf{\int_M{|\nabla u_i|^2dV}}},$$

\nd and
$$\ddd\int_M{au^2dV} \ = \ \ddd\lim_{i \to \infty}{\int_M{a(u_i)^2dV}}.$$

Thus,

\begin{eqnarray*}
I(u) & = &  \ddd\frac{1}{2}\int_M{|\nabla u|^2dV} \ + \ \ddd\frac{1}{2}\int_M{au^2dV} \ - \ \bigl<\Psi, u\bigr> \\
     & \leq &  \ddd\frac{1}{2}\lim_{i \to \infty}{\inf{\int_M{|\nabla u_i|^2dV}}} \ + \ \ddd\frac{1}{2}\lim_{i \to 
\infty}{\int_M{au_i^2dV}} \ - \ \ddd\lim_{i \to \infty}{\bigl<\Psi, u_i\bigr>} \\
  & = & \ddd\lim_{i \to \infty}{\inf{I(u_i)}} 
\end{eqnarray*}

what prove the Claim 1. \fimaf
\vsp

As \ $H^2_1$ \ is a reflexive space, (see \cite{kavian}) there is \ $u \in H^2_1$, \ such that 
$$I(u) \ = \ \ddd\min_{v \in H^2_1}{I(v)}$$ 

\nd and, consequently, \ $u$ \ is a solution of the equation  \ (\ref{4.1}).

{{\bf Claim 2} \it
The uniqueness of solution  for the equation  \ (\ref{4.1}) holds.
}

Let us suppose that \ $u_1$ \ and \ $u_2$ \ are solutions of the  equation  \ (\ref{4.1}). Then,
 \ $\forall \ \varphi \ \in \ H^2_1$, \ we have that 

\begin{equation}
 \label{4.2}
\ddd\int_M{\nabla u_1\nabla \varphi dV} \ + \ \ddd\int_M{au_1\varphi dV} \ = \ \bigl<\Psi, \varphi \bigr>,
\end{equation}

\begin{equation}
\label{4.3}
\ddd\int_M{\nabla u_2\nabla \varphi dV} \ + \ \ddd\int_M{au_2\varphi dV} \ = \ \bigl<\Psi, \varphi \bigr>.
\end{equation}

\nd Taking \ $\varphi \ = \ u_1 \ - \ u_2$ \ and considering the difference  (\ref{4.2}) \ and (\ref{4.3}) \ we get
$$\ddd\int_M{|\nabla (u_1 - u_2)|^2dV} \ + \ \ddd\int_M{a(u_1 - u_2)^2dV} \ = 0.$$

\nd Thus, as \ $a \ > \ 0$, then \ $u_1 \ = \ u_2$ \ in \ $H^2_1$. What give us the proof of Claim 2. \fimaf

{{\bf Claim 3} \it
The operator $T$ \ is  continuous.
}

Indeed, let \ $\{\Psi_i\} \ \ \mbox{and} \ \Psi \ \in \ (H^2_1)^{'}$  such that 
$$\Psi_i \ \longrightarrow \ \Psi \ \mbox{in} \ (H^2_1)^{'}.$$ 

\nd Taking 

$$u_i \ = \ T(\Psi_i) \ \ \mbox{and} \ u \ = \ T(\Psi),$$ 

\nd obtain \ $\forall \ \varphi \ \in \ H^2_1,$
\begin{equation}
 \label{4.4}
\ddd\int_M{\nabla u_i\nabla \varphi dV} \ + \ \ddd\int_M{au_i\varphi dV} \ = \ \bigl<\Psi_i, \varphi \bigr>,
\end{equation}
\begin{equation}
\label {4.5}\ddd\int_M{\nabla u\nabla \varphi dV} \ + \ \ddd\int_M{au\varphi dV} \ = \ \bigl<\Psi, \varphi \bigr>.
 \end{equation}

Substituting \ $\varphi_i \ = \ u_i \ - \ u \ \in  \ H^2_1$ \ in \ (\ref{4.4}) and \ (\ref{4.5}) we get that
\begin{equation}
 \label{4.6}
\ddd\int_M{\nabla u_i(\nabla u_i - \nabla u) dV} \ + \ \ddd\int_M{au_i(u_i - u)dV} \ = \ \bigl<\Psi_i, u_i - u\bigr>,
\end{equation}
    \begin{equation}
     \label{4.7}
\ddd\int_M{\nabla u(\nabla u_i - \nabla u) dV} \ + \ \ddd\int_M{au(u_i - u)dV} \ = \ \bigl<\Psi, u_i - u\bigr>.
    \end{equation}

Considering \ (\ref{4.6}) \ - \ (\ref{4.7})  we obtain  

\begin{eqnarray*}
\ddd\int_M{|\nabla (u_i - u)|^2 dV} \  + \  \ddd\int_M{a(u_i - u)^2dV} & = & \bigl<\Psi_i - \Psi, u_i - u\bigr> \\
 &  \leq & \|\Psi_i - \Psi\|_{(H^2_1)^{'}}\| u_i - u\|_{H^2_1}.   
\end{eqnarray*}

As 
$$\ddd\int_M{|\nabla (u_i - u)|^2 dV} \ + \ \ddd\int_M{a(u_i - u)^2dV} \ \geq \ C\| u_i - u\|^2_{H^2_1},$$

\nd where  
$$C \ = \ \min{\{1,\min a\}} \ > \ 0,$$

\nd we obtain
$$\| u_i - u\|_{H^2_1} \ \leq \ \ddd\frac{1}{C}\|\Psi_i - \Psi\|_{(H^2_1)^{'}} \ \longrightarrow \ 0.$$

Therefore, 
$$u_i \ \longrightarrow \ u \ \mbox{in} \ H^2_1,$$

\nd namely, \ $T$ \ is continuous. Proving the Claim 3. \fimaf

{{\bf Claim 4} \it The operator
$T$ \ is non-decreasing.
}

Let \ $\Psi_1, \ \Psi_2 \ \in \ (H^2_1)^{'}$ \ such that \ $\Psi_1 \ \leq \ \Psi_2$ \ in the sense of  that 
$$\bigl<\Psi_1, \varphi \bigr> \ \ \leq \ \ \bigl<\Psi_2, \varphi \bigr>, \ \ \forall \ \varphi \ \in \ H^2_1, \ \varphi
 \ \geq \ 0.$$

Taking \ $u_1 \ = \ T(\Psi_1) \ \ \mbox{and} \ u_2 \ = \ T(\Psi_2) \ \in \ H^2_1$, we want to show that \ $u_1 \ \leq
 \ u_2$ \ a.e. in \ $M$. \ Indeed, as for all \ $\varphi \ \in \ H^2_1, \ \varphi \ \geq \ 0,$

\begin{eqnarray*}
\ddd\int_M{\nabla u_1\nabla \varphi dV} \ + \ \ddd\int_M{au_1\varphi dV} & = & \bigl<\Psi_1, \varphi \bigr> \\
 & \leq & \bigl<\Psi_2, \varphi \bigr> \\ 
 & = & \ddd\int_M{\nabla u_2\nabla \varphi dV} \ + \ \ddd\int_M{au_2\varphi dV}.
\end{eqnarray*} 

It follows by the Weak Comparison Principle that 
\vspace{-0.5mm}
$$u_1 \ \leq \ u_2 \ \mbox{a.e. in} \ M.$$ \fimaf 
\vsp

The proof of Lemma 1 follows immediately using Claims 1-4. \fim
\vsp

\section{The Proof of the Theorem 1}

Now we will proof the Main Result.

Consider the ``interval" 
$$[\underline{u},\overline{u}] \ = \ \{v \ \in \ H^2_1; \ \underline{u}(x) \ \leq \ v(x) \ \leq \ \overline{u}(x) \ 
\mbox{a.e. in} \ M \}$$

\nd with the topology of the  a.e. convergence, consider 
$$S \ : \ [\underline{u},\overline{u}] \ \longrightarrow \ (H^2_1)^{'}$$ 

\nd by 
$$\bigl<S(v),\varphi \bigr> \ = \ \ddd\int_M{f(x)F(v) \varphi dV} \ + \ \ddd\int_M{h(x)H(v)\varphi dV}, \ v \ \in \ 
[\underline{u},\overline{u}], \ \varphi \ \in \ H^2_1.$$

{{\bf Claim 5} \it
$S$ \ is a continuous, non-decreasing and bounded operator.
}
 
\nd Proof of Claim 5:
 
\nd $(i)$ \ $S$ \ is bounded.
\vsp

Indeed, if \ $v \ \in \ [\underline{u},\overline{u}] \ \mbox{and} \  \varphi \ \in \ H^2_1$, we have that

\begin{eqnarray*}
|\bigl<S(v),\varphi \bigr>| & \leq & \ddd\int_M{f(x)F(\overline{u})| \varphi| dV} \ + \ \ddd\int_M{h(x)H(\overline{u})|
\varphi| dV} \\
& \leq & \|fF(\overline{u})\|_2\|\varphi \|_2 \ + \  \|hH(\overline{u})\|_2\|\varphi \|_2 \\
& = & \left( \|fF(\overline{u})\|_2 \ + \  \|hH(\overline{u})\|_2\right) \|\varphi \|_2
\end{eqnarray*}

\nd where in the last inequality we use the H$\ddot{o}$lder's inequality.
\vsp

As \ $H^2_1 \ \hookrightarrow \ L^2$, \ we obtain that
$$|\bigl<S(v),\varphi \bigr>| \ \leq \ \left( \|fF(\overline{u})\|_2 \ + \  \|hH(\overline{u})\|_2\right)C
 \|\varphi \|_{H^2_1} \ = \ A \|\varphi \|_{H^2_1},$$

\nd where \ $A \ = \ C\left( \|fF(\overline{u})\|_2 \ + \  \|hH(\overline{u})\|_2\right) \ > \ 0$.
\vsp

Hence, \ $S$ \ is bounded in \ $[\underline{u},\overline{u}]$.
\vsp

\nd $(ii)$ \ $S$ \ is non-decreasing.
\vsp

Indeed, if \ $u_1, \ u_2 \ \in \ [\underline{u},\overline{u}]$ \ are such that \ $u_1 \ \leq \ u_2$ \ a.e., it follows that,
 by the fact that \ $F \ \mbox{and} \ H$ \ are non-decreasing and by \ $(\alpha_2)$, that for all \ $\varphi \ \in \ H^2_1,
 \ \varphi \ \geq \ 0$, 

\begin{eqnarray*}
\bigl<S(u_1),\varphi \bigr> & = & \ddd\int_M{f(x)F(u_1) \varphi dV} \ + \ \ddd\int_M{h(x)H(u_1)\varphi dV} \\
& \leq & \ddd\int_M{f(x)F(u_2) \varphi dV} \ + \ \ddd\int_M{h(x)H(u_2)\varphi dV} \\
& = & \bigl<S(u_2),\varphi \bigr>.
\end{eqnarray*}

\nd $(iii)$ \ $S$ \ is continuous.
\vsp

Let \ $(v_i) \ \mbox{and} \ v \ \in \ [\underline{u},\overline{u}]$ \ such that \ $v_i \ \longrightarrow \ v$ \ a.e. in \ $M$. 
\vsp

We observe that \ $fF(v_i) \ + \ hH(v_i) \ \longrightarrow \ fF(v) \ + \ hH(v)$ \ a.e. in \ $M$ \ and 

\begin{equation}
 \label{4.8}
|fF(v_i) \ + \ hH(v_i)|^2 \ \leq \ |fF(\overline{u}) \ + \ hH(\overline{u})|^2.
\end{equation}

\vsp

On the other hand, for \ $\varphi \ \in \ H^2_1$

\begin{eqnarray*}
|\bigl<S(v_i) - S(v),\varphi \bigr>| & \leq & \ddd\int_M{|(fF(v_i) + hH(v_i)) - (fF(v) + hH(v))|| \varphi| dV} \\
& \leq &  \|(fF(v_i) + hH(v_i)) - (fF(v) + hH(v))\|_2\|\varphi \|_2 \\
& \leq & C\|(fF(v_i) + hH(v_i)) - (fF(v) + hH(v))\|_2 \|\varphi \|_{H^2_1}
\end{eqnarray*}

\nd where  in the two last inequalities we used the H$\ddot{o}$lder's inequality and the Sobolev's embedded \ $H^2_1 \ \hookrightarrow 
\ L^2$, \ respectively.
\vsp

Then, by equation \ (\ref{4.8})  we can apply the Theorem of Dominated Convergence of Lebesgue to conclude that 
$$\|S(v_i) - S(v) \|_{(H^2_1)^{'}} \ = \ o(1).$$

Therefore, \ $S$ \ is continuous in \ $[\underline{u},\overline{u}]$ \ what proves the Claim 5. \fimaf
\vsp

Consider, now, \ $J \ : \ [\underline{u},\overline{u}] \ \longrightarrow \ H^2_1$ \ given by \ $J \ = \ T\circ S$. 
\ Namely, for all \ $v \ \in \ [\underline{u},\overline{u}], \ u \ = \ J(v)$ \ is the unique solution of equation
$$\Delta u \ + \ au \ = \ fF(v) \ + \ hH(v).$$

Taking 
$$u_1 \ = \ J(\underline{u}), \ u^1 \ = \ J(\overline{u}), \ u_{i+1} \ = \ J(u_i) \ \mbox{and} \  u^{i+1} \ = \ J(u^i), \ \ i
 \ \geq \ 1,$$ 

\nd we obtain 
$$\Delta u_1 \ + \ au_1 \ = \ fF(\underline{u}) \ + \ hH(\underline{u}) \ \geq \ \Delta \underline{u} \ + \ a\underline{u}$$

and by the Weak Comparison Principle \ $u_1 \ \geq \ \underline{u}$ \ a.e.
\vsp

Analogously, 
$$\Delta u_2 \ + \ au_2  \ = \ fF(u_1) \ + \ hH(u_1) \ \geq \  fF(\underline{u}) \ + \ hH(\underline{u}) \ = \ \Delta u_1
 \ + \ au_1$$

and, thus, \ $u_2 \ \geq \ u_1 $ \ a.e.
\vsp

Considering this process, successively, we obtain that 
$$\underline{u} \ \leq u_1 \ \leq \ u_2 \ \leq \ \ldots \ \leq \ u_i \ \leq \ \ldots $$

With the same argument we conclude that
$$\ldots \ \leq u^i \ \leq \  \ldots \ \leq \ u^2 \ \leq \ u^1 \ \leq \ \overline{u}. $$

On the other hand,
$$\Delta u^1 \ + \ au^1  \ = \ fF(\overline{u}) \ + \ hH(\overline{u}) \ \geq \  fF(\underline{u}) \ + \ hH(\underline{u}) \
= \ \Delta u_1 \ + \ au_1$$

\nd what give us \ $u^1 \ \geq \ u_1 \ \geq \ \underline{u}$.
\vsp

Analogously,
$$\Delta u^2 \ + \ au^2  \ = \ fF(u^1) \ + \ hH(u^1) \ \geq \  fF(u_1) \ + \ hH(u_1) \ = \ \Delta u_2 \ + \ au_2$$

\nd then \ $u^2 \ \geq \ u_2$. 
\vsp

Considering this process, successively, we obtain that 
$$\underline{u} \ \leq \ u_1 \ \leq \ u_2 \ \leq \ \ldots \ \leq \ u_i \ \leq \ \ldots  \  \leq u^i \ \leq \  \ldots 
\ \leq \ u^2 \ \leq \ u^1 \ \leq \ \overline{u}. $$

Then, \ $u_i \ \longrightarrow \ u_{\ast} \ \ \mbox{and} \ \ u^i \ \longrightarrow \ u^{\ast}$ \ a.e. when \ $i 
\ \rightarrow \ \infty$ \ and \ $u_{\ast}, \ u^{\ast} \ \in \ [\underline{u},\overline{u}] \ \ \mbox{ with} \ u_{\ast} 
\ \leq  \ u^{\ast}$ \ a.e.
\vsp

As \ $u_{i+1} \ = \ J(u_i) \ \longrightarrow \ J(u_{\ast}) \ \ \mbox{and} \ \ u^{i+1} \ = \ J(u^i) \ \longrightarrow \ 
J(u^{\ast}) \ \mbox{ when} \ i \ \longrightarrow \ \infty$, \ by continuity of \ $J$, \ we conclude that \ $u_{\ast}, 
\ u^{\ast} \ \in \ H^2_1$ \ with \ $J(u_{\ast}) \ = \ u_{\ast} \ \ \mbox{ and} \ \ J(u^{\ast}) \ = \ u^{\ast}$, \ namely, \ 
$\underline{u} \ \leq \ u_{\ast} \ \leq \ u^{\ast}  \leq \ \overline{u} $ \ are weak solutions of equation
$$ \Delta u \ + \ au \ = \ fF(u) \ + \ hH(u).$$
By using \ $(\alpha_1)$ \ we can use a Regularity Theorem to  show that \ $u_{\ast}, \ u^{\ast} \ \in
 \ C^{\infty}(M)$ \ and by the Strong Maximum Principle \ $u_{\ast}, \ u^{\ast} \ > \ 0$ \ in \ $M$.  \fim

{INSTITUTO DE CI\^ENCIAS EXATAS E DA TERRA\\
 CAMPUS UNIVERSIT\'ARIO DO ARAGUAIA\\UNIVERSIDADE FEDERAL DE MATO GROSSO
 \\  78698-000 Pontal do Araguaia,  MT, Brazil}
 
email ={carlosro@ufmt.br}

INSTITUTO DE MATEMÁTICA E ESTATÍSTICA\\
UNIVERSIDADE FEDERAL DE GOIÁS - UFG\\
74690-900, Goi\^ania, GO, Brazil

email= {msouza@ufg.br}

\end{document}